\documentclass[11pt,a4paper,twoside,english]{article}
\usepackage{color}
\usepackage{latexsym,amsmath,amssymb} 
\usepackage{babel}
\usepackage{graphicx} 
\title{\bf \Large Equivalence between
limit theorems for lattice group-valued $k$-triangular set functions}
\author{A. Boccuto \thanks{
Dipartimento di Matematica e Informatica, University of Perugia,
via Vanvitelli 1, I-06123 Perugia, Italy; e-mail:
antonio.boccuto@unipg.it, boccuto@yahoo.it (Corresponding author)}
\and  X. Dimitriou\thanks{ Department of Mathematics, University of Athens,
Panepistimiopolis, Athens 15784, Greece; $\qquad$ $\qquad$
$\qquad$ $\qquad$ email: xenofon11@gmail.com,
dxenof@math.uoa.gr  \newline 
 \textit{2010 A. M. S. Subject Classifications:}
Primary: 26E50, 28A12, 28A33, 28B10, 28B15,
40A35, 46G10, 54A20, 54A40. 
Secondary: 06F15, 06F20, 06F30, 22A10, 
28A05, 40G15, 46G12, 54H11, 54H12. \newline {\it Key words}:
lattice group, $(D)$-convergence,  
$(O)$-convergence, $k$-subadditive set function, 
$k$-triangular set function, continuous set function,
$(s)$-bounded set function, Fremlin lemma, limit theorem,
Brooks-Jewett theorem, Vitali-Hahn-Saks theorem, Nikod\'{y}m 
theorem, Drewnowski theorem, Schur theorem.}}
\date{}
\begin{document}
\maketitle 
\newtheorem{theorem}{Theorem}[section]
\newtheorem{corollary}[theorem]{Corollary}
\newtheorem{lemma}[theorem]{Lemma} 
\newtheorem{proposition}[theorem]{Proposition}
\newtheorem{definition}[theorem]{Definition}
\newtheorem{definitions}[theorem]{Definitions}
\newtheorem{remark}[theorem]{Remark}
\newtheorem{remarks}[theorem]{Remarks} 
\newtheorem{example}[theorem]{Example}
\newcommand{\bvi}{\bigvee_{t=1}^{\infty} a_{t, \varphi(t)}}
\newcommand{\bviastar}{\bigvee_{t
=1}^{\infty} a^*_{t, \varphi(t)}}
\newcommand{\bvib}{\bigvee_{t=1}^{\infty} b_{t, \varphi(t)}}
\newcommand{\bvialpha}{\bigvee_{t=1}^{\infty} \alpha_{t, \varphi(t)}}
\newcommand{\bvic}{\bigvee_{t=1}^{\infty} c_{t, \varphi(t)}}
\newcommand{\bvid}{\bigvee_{t=1}^{\infty} d_{t, \varphi(t)}}
\newcommand{\bviA}{\bigvee_{t=1}^{\infty} A_{t, \varphi(t)}}
\newcommand{\bviB}{\bigvee_{t=1}^{\infty} B_{t, \varphi(t)}}
\newcommand{\bviC}{\bigvee_{t=1}^{\infty} C_{t, \varphi(t)}}
\newcommand{\bviD}{\bigvee_{t=1}^{\infty} D_{t, \varphi(t)}}
\newcommand{\bvie}{\bigvee_{t=1}^{\infty} e_{t, \varphi(t)}}
\newcommand{\bviz}{\bigvee_{t=1}^{\infty} z_{t, \varphi(t)}}
\begin{abstract} 
We investigate some main properties of 
lattice group-valued $k$-triangular set functions
and prove some Brooks-Jewett,
Nikod\'{y}m, Vitali-Hahn-Saks and Schur-type theorems
and their equivalence. A Drewnowski-type theorem on 
existence of continuous restrictions of $(s)$-bounded set 
functions is given. 
Furthermore we pose some open problems.
\end{abstract}
\section{Introduction}
In this paper we deal with
$k$-triangular lattice group-valued set 
functions. We continue the investigation started in 
\cite{bdsubadditive, bdgreece2016, bdgreece2017}, 
where some limit theorems were proved for
$k$-subadditive, positive and monotone set functions.
In particular we treat $(s)$-boundedness and 
continuity from above at $\emptyset$. 
Among the related literature, see for instance 
\cite{BSROMA, 
guariglia1990, 
papvhs, PAPBOOK, SAEKI, ventriglianapoli,
WANGKLIR}.
Some examples of $k$-triangular set functions are the 
so-called ``$M$-measures'', namely 
monotone set functions $m$ with $m(\emptyset)=0$, 
continuous from above and from below and compatible
with respect to finite suprema and infima, which have 
several applications, for example to intuitionistic fuzzy 
sets and observables (see also 
\cite{BAN, 
RN}). A class of $1$-triangular set functions which are not 
necessarily monotone
is that of the Saeki measuroids (see also \cite{SAEKI}).
We investigate some main properties 
and prove some Brooks-Jewett, Nikod\'{y}m, 
Vitali-Hahn-Saks and Schur-type theorem and their 
equivalence,
by means of sliding hump-type techniques. Here 
the notions of (uniform) $(s)$-boundedness,
continuity from above at $\emptyset$ and continuity with
respect to a suitable Fr\'{e}chet-Nikod\'{y}m topology
$\tau$, as well as the concept of 
pointwise convergence of the involved set functions, are 
intended with respect to a single regulator. We deal with 
$(D)$-convergence
and use the Fremlin theorem, a powerful tool which allows to
reduce a sequence of regulators to a single $(D)$-sequence.
We prove a Drewnowski-type theorem 
about the existence of continuous 
restrictions of $(s)$-bounded $k$-triangular set functions,
extending earlier results proved in \cite{bcdefin, bdequiv, 
DRE2}. Note that in our setting, differently from  
\cite{bdequiv}, it is possible to use a 
Drewnowski-type technique just because the involved 
$(D)$-convergence is always meant with respect to a single
$(D)$-sequence, and so it is possible to prove that the
Nikod\'{y}m convergence theorem implies the Brooks-Jewett theorem,
analogously as in \cite{DRE2}, just by choosing a disjoint set sequence
and considering the behavior of suitable subsequences. 
We use a sliding-hump type technique 
(see also \cite{CANDELOROPALERMO, CANDELOROLETTA}) 
Furthermore, observe that here $(s)$-boundedness of the
limit set function
is proved as thesis, while in \cite[Theorem 1]{papvhs} 
and \cite[Theorem 11.10]{PAPBOOK} it is assumed as hypothesis.
Some other versions of limit theorems for lattice group-valued set 
functions in the recent literature are proved, for instance, in
\cite{KADETS, 
bdequiv, 
bdpmodena, 
bdpschurfilters, cs}
(see also \cite{bdbook, bdLAP} for an overview).
Finally we pose some open problems.
\section{Preliminaries}
We begin with recalling the following main properties of lattice
groups (see also \cite{bdbook, bdLAP, BRV}).
\begin{definitions} \rm
(a) A Dedekind complete lattice group $R$ is said to be \textit{super
Dedekind complete} iff for every nonempty set
$A \subset R$, bounded from above, there is a 
countable subset $A^{\prime}$, with $\displaystyle{
\bigvee A^{\prime}=\bigvee A}$.

(b) A sequence $(\sigma_p)_p$ in $R$ is called an 
\textit{$(O)$-sequence} iff it is decreasing and 
$\displaystyle{\bigwedge_p \, \sigma_p=0}$. 

(c) A bounded double sequence $(a_{t,l})_{t,l}$ in
$R$ is a \textit{$(D)$-sequence} or a
\textit{regulator} iff $(a_{t,l})_l$ is an
$(O)$-sequence for any $t \in \mathbb{N}$.

(d) A lattice group $R$ is \textit{weakly
$\sigma$-distributive} iff
$\displaystyle{\bigwedge_{\varphi \in
{\mathbb{N}}^{\mathbb{N}} } \Bigl(
\bigvee_{t=1}^{\infty}
a_{t,\varphi(t)}\Bigr)=0}$
for every $(D)$-sequence $(a_{t,l})_{t,l}$ in $R$.

(e) A sequence $(x_n)_{n}$ in $R$ is said to be \em order convergent 
\rm (or $(O)$\em -convergent \rm ) to $x$ iff there exists an
$(O)$-sequence $(\sigma_p)_p$ in $R$ such that for every 
$p\in\mathbb{N}$ there is a positive integer $n_0$ with
$|x_n - x| \leq \sigma_p$ for
each $n \geq n_0$, and in this case we write
$\displaystyle{(O)\lim_n\,x_n = x}$.

(f) A sequence $(x_n)_n$ in $R$ is $(D)$\em -convergent \rm to $x$ iff 
there is a $(D)$-sequence $(a_{t,l})_{t,l}$ in $R$ such that for every
$\varphi\in\mathbb{N}^{\mathbb{N}}$ there is $n_0
\in \mathbb{N}$ with $\displaystyle{
|x_n - x| \leq\bigvee_{t=1}^{\infty} a_{t,\varphi(t)} 
}$ whenever $n \geq n_0$, and we write 
$\displaystyle{(D)\lim_n\,x_n = x}$.

(g) We call \textit{sum} of a series 
$\displaystyle{\sum_{n=1}^{\infty} x_n}$ in $R$ the limit 
$\displaystyle{(O)\lim_n \sum_{r=1}^n x_r}$, if it
exists in $R$.
\end{definitions}
\begin{remark}
\rm Observe that in every Dedekind complete lattice group $R$ any
$(O)$-convergent sequence is $(D)$-convergent too, while 
the converse is true if and only if $R$ is weakly
$\sigma$-distributive (see also \cite{bdLAP, BRV}).
\end{remark}
We now recall the Fremlin lemma, which has a fundamental 
importance 
in the setting of $(D)$-convergence, because it allows to 
replace a sequence of regulators with a single $(D)$-sequence.
\begin{lemma}\label{fremlin} {\rm (%
see also \cite[Lemma 1C]{FREMLIN},
\cite[Theorem 3.2.3]{RN})} \em
Let $R$ be any Dedekind complete 
lattice group
and $(a_{t,l}^{(n)})_{t,l}$, $n \in
\mathbb{N}$, be a sequence of regulators in $R$.
Then for every $u\in R$, $u \geq 0$ there is a
$(D)$-sequence $(a_{t,l})_{t,l}$ in $R$ with
\begin{eqnarray*}
\displaystyle{u\wedge\Bigl(\sum_{n=1}^{q}
\Bigl(\bigvee_{t=1}^{\infty} a^{(n)}_{t,
\varphi(t+n)}\Bigr)\Bigr)\leq\bigvee_{t=1}^{\infty}
a_{t,\varphi(t)}}\quad \text{for every  }q \in
\mathbb{N}\text{  and   } \varphi \in
\mathbb{N}^{\mathbb{N}}.
\end{eqnarray*}
\end{lemma}
The following result links $(O)$-sequences and regulators.
\begin{theorem}\label{bonitra} \rm (see also
\cite[Theorem 3.4]{bonitra}) \em
Given any De\-de\-kind complete lattice group $R$
and any $(O)$-sequence $(\sigma_l)_l$ in $R$,
the double sequence defined by
$a_{t,l}:=\sigma_l$, $t$, $l\in\mathbb{N}$, is a
$(D)$-sequence, and for every
$\varphi\in{\mathbb{N}}^{\mathbb{N}}$,
\begin{eqnarray}\label{bonitra1}
\sigma_l\leq
\bigvee_{t=1}^{\infty}\,a_{t,\varphi(t)},
\end{eqnarray} where $l=\varphi(1)$.
Conversely, if $R$ is super Dedekind complete and
weakly $\sigma$-distributive, then for
any regulator $(a_{t,l})_{t,l}$ in $R$
there is an $(O)$-sequence $(\sigma_p)_p$ such
that for each $p \in \mathbb{N}$ there is
$\varphi_p \in {\mathbb{N}}^{\mathbb{N}}$
with \begin{eqnarray}\label{bonitra2}
\bigvee_{t=1}^{\infty}\, a_{t,\varphi_p(t)}\leq \sigma_p .
 \end{eqnarray}
\end{theorem} 

We now deal with some fundamental properties of
non-additive
lattice group-valued set functions (see also 
\cite{bdLAP, PAPBOOK, 
WANGKLIR} and the bibliographies therein). 

From now on, when it is not indicated elsewhere 
explicitly, $R$ denotes a Dedekind complete and
weakly $\sigma$-distributive
lattice group, $R^+$ is the positive cone of $R$,
$G$ is an infinite set, ${\mathcal L}$ is an algebra
of subsets of $G$, $m:{\mathcal L} \to R$ is a 
positive bounded set function
and $k$ is a fixed positive integer.
\begin{definitions} \rm 
(a) We say that $m$ is 
\textit{$k$-subadditive} on ${\mathcal L}$ iff $m(\emptyset)=0$ and
\begin{eqnarray}\label{ksubadditivity}
m(A\cup B) \leq m(A) + k\, m(B) \quad
\text{  whenever  } A, B\in{\mathcal L}, \, A \cap B=\emptyset ;
\end{eqnarray}
\textit{$k$-triangular} on ${\mathcal L}$, iff
$m$ is $k$-subadditive and  
\begin{eqnarray}\label{ktriangularity}
m(A) - k\, m(B) \leq m(A \cup B) \quad \text{ whenever  }
A , B \in {\mathcal L}, \, \, A \cap B=\emptyset. 
\end{eqnarray}

(b) We call 
\textit{semivariation} of $m$, shortly $v(m)$, the set function
defined by 
\begin{eqnarray}\label{semivariation}
v(m) (A):=\bigvee \{m(B): B \in  {\mathcal L}, \, \, B \subset A \},
\quad A \in {\mathcal L}. 
\end{eqnarray}
\end{definitions}
The following results will be useful in the sequel.
\begin{proposition}\label{semivinherit}
If $m:{\mathcal L} \to R$ is $k$-subadditive, then 
$v(m)$ is $k$-triangular.
\end{proposition}
{\bf Proof:} Pick any two disjoint sets $A$, $B \in
{\mathcal L}$. Fix arbitrarily $C \subset A \cup B$ and set 
$C_1:=A \cap C$, $C_2:=B \cap C$. By $k$-subadditivity of 
$m$ we get
\begin{eqnarray*}\label{k1}
m(C) \leq m(C_1) + k\, m(C_2),
\end{eqnarray*} and  hence
\begin{eqnarray}\label{k3}
m(C) \leq v(m)(A) + k\, v(m)(B)
\end{eqnarray} By arbitrariness of $C$, from
(\ref{k3}) we obtain $k$-subadditivity of $v(m)$. 
Relation (\ref{ktriangularity}) 
follows from monotonicity of $v(m)$,
taking into account that $v(m)(\emptyset)=0$. $\quad \Box$
\begin{proposition}\label{kfinitetriangularity} Let 
$m:{\mathcal L} \to R$ be a $k$-triangular set function. Then 
for every $n\in\mathbb{N}$, $n\geq 2$, and for every pairwise
disjoint sets $E_1$, $E_2,
\ldots, E_n \in {\mathcal L}$ it is
\begin{eqnarray}\label{kntriangularity}
m(E_1) - k \sum_{q=2}^n m(E_q) \leq 
m\Bigl( \bigcup_{q=1}^n E_q \Bigr) \leq m(E_1)+ k \sum_{q=2}^n
m(E_q).
\end{eqnarray} In particular we have
\begin{eqnarray}\label{kntriangularitybis}
m(E_1) \leq m\Bigl( \bigcup_{q=1}^n E_q \Bigr) + k \sum_{q=2}^n
m(E_q).
\end{eqnarray}
\end{proposition}
{\bf Proof:} Let us proceed by induction, and assume that 
(\ref{kntriangularity}) holds for $n-1$. We get
\begin{eqnarray*}\label{finalinitial}
m(E_1)-k \sum_{q=2}^n m(E_q)&=&m(E_1)-k 
\sum_{q=2}^{n-1}m(E_q)-
k \, m(E_n) \leq \nonumber
\\&\leq& m\Bigl( \bigcup_{q=1}^{n-1}E_q \Bigr) - k \, m(E_n) \leq
m\Bigl( \bigcup_{q=1}^n E_q \Bigr) \leq m\Bigl( 
\bigcup_{q=1}^{n-1}E_q \Bigr) + k \, m(E_n)  \nonumber
\leq\\&\leq&m(E_1)+k \sum_{q=2}^{n-1}m(E_q)+
k \, m(E_n)=m(E_1)+k \sum_{q=2}^n m(E_q),  \nonumber
\end{eqnarray*} obtaining the assertion. $\quad \Box$
\begin{definitions} \rm (a) 
A lattice ${\mathcal E}$ of
subsets of $G$ is said to 
\textit{satisfy property $(E)$} iff every disjoint 
sequence $(C_h)_h$ in ${\mathcal E}$
 has a subsequence $(C_{h_r})_r$,
such that ${\mathcal E}$ contains 
the $\sigma$-algebra 
generated by the sets $C_{h_r}$, $r\in\mathbb{N}$, in the set
$\displaystyle{ \bigcup_{r=1}^{\infty} C_{h_r}}$ (see also
\cite{SCHACHERMAYER}).

From now on we assume that ${\mathcal L} \subset 
{\mathcal P}(G)$ is an algebra, satisfying property $(E)$.

(b) A set function $m: {\mathcal L} \to R^+$ is
\textit{$(s)$-bounded} on ${\mathcal L}$ iff there exists a 
$(D)$-sequence $(a_{t,l})_{t,l}$ such that, 
for every disjoint sequence $(C_h)_h$ in $
{\mathcal L}$,
$\displaystyle{(D)\lim_h m(C_h)=0}$
with respect to $(a_{t,l})_{t,l}$.

(c) We say that the set functions $m_j: {\mathcal L} \to R^+$, 
$j \in \mathbb{N}$, are \textit{uniformly
$(s)$-bounded} on ${\mathcal L}$ iff there exists a 
$(D)$-sequence $(a_{t,l})_{t,l}$ such that, 
for every disjoint sequence $(C_h)_h$ in $
{\mathcal L}$,
$$\displaystyle{(D)\lim_h \Bigl(\bigvee_j m_j(C_h)\Bigr)=0}$$
with respect to $(a_{t,l})_{t,l}$.

(d) We say that a set function 
$m:{\mathcal L} \to R^+$ 
is \textit{continuous from above at
$\emptyset$} iff there is a $(D)$-sequence $(a_{t,l})_{t,l}$ with
$\displaystyle{(D)\lim_n m(H_n)=0}$ with respect to 
$(a_{t,l})_{t,l}$ 
whenever $(H_n)_n$ is a decreasing sequence in ${\mathcal L}$ with
$\displaystyle{\bigcap_{n=1}^{\infty}H_n=\emptyset}$.

(e) The set functions $m_j:{\mathcal L} \to R^+$, $j\in\mathbb{N}$, are said 
to be \textit{uniformly continuous from above at $\emptyset$} iff
there is a $(D)$-sequence $(a_{t,l})_{t,l}$ with
$$(D)\lim_n \Bigl( \bigvee_j v(m_j)(H_n) \Bigr)=0$$ 
with respect to $(a_{t,l})_{t,l}$
for each decreasing sequence
$(H_n)_n$ of elements of ${\mathcal L}$ with
$\displaystyle{\bigcap_{n=1}^{\infty}H_n=\emptyset}$.

(f) We say that the set functions $m_j: {\mathcal L}\to R$, $j\in
\mathbb{N}$, are \textit{equibounded} on ${\mathcal L}$
iff there is an element $u\in 
R$ with $|m_j(A)| \leq u$ for all $j\in \mathbb{N}$ and $A \in 
{\mathcal L}$. 

\end{definitions}
\begin{remark}\label{29}
 \rm Observe that continuity from above at 
$\emptyset$ of a $k$-triangular set function with respect to a 
regulator $(a_{t,l})_{t,l}$ implies its
$(s)$-boundedness with respect to the $(D)$-sequence 
$( (k+1) a_{t,l})_{t,l}$. 
Otherwise there exist a disjoint
sequence $(C_h)_h$ in ${\mathcal L}$ and
an element $\varphi\in{\mathbb{N}}^{\mathbb{N}}$ with
\begin{eqnarray}\label{1515}
m(C_h) \not \leq (k+1) \bigvee_{t=1}^{\infty} a_{t,\varphi(t)}
\quad \text{  for every   } h\in \mathbb{N}.
\end{eqnarray}
Since ${\mathcal L}$ satisfies property $(E)$, 
possibly passing to a suitable subsequence of $(C_h)_h$,
without loss of generality we can suppose that 
$$\bigcup_{h\in P} C_h \in {\mathcal L} \quad \text{ for  all  }
P \subset \mathbb{N}.$$ As $m$ is $k$-triangular and
continuous from above at $\emptyset$ with respect to the 
$(D)$-sequence $(a_{t,l})_{t,l}$, then for each $\varphi\in
{\mathbb{N}}^{\mathbb{N}}$ there exists a positive integer $h_0$
with \begin{eqnarray*}
m(C_h) \leq m \Bigl( \bigcup_{i=h}^{\infty} C_i \Bigr) + k \, 
m \Bigr( \bigcup_{i=h+1}^{\infty} C_i \Bigl) \leq (k+1)
\bigvee_{t=1}^{\infty}a_{t,\varphi(t)} \text{  for  any  } h \geq h_0,   
\end{eqnarray*} getting a contradiction with (\ref{1515}).
Analogously it is possible to check that uniform continuity from 
above at $\emptyset$ implies uniform $(s)$-boundedness.

On the other hand, in general the converse implication is not true.
For example observe that a real-valued measure $m$ defined
on a $\sigma$-algebra $\Sigma \subset {\mathcal P} (G)$   
is countably additive (resp. $(s)$-bounded) if and only if the
set function $A\mapsto |m(A)|$, $A\in \Sigma$, which is by 
construction $1$-triangular on $\Sigma$, is continuous from above at 
$\emptyset$ (resp. $(s)$-bounded) (see also 
\cite[Example (a)]{SAEKI}).
\end{remark}
We now prove that the semivariation inherits both 
$(s)$-boundedness and continuity from above at $\emptyset$.
\begin{proposition}\label{3434}
Let $m:{\mathcal L} \to R^+$ be an $(s)$-bounded set function.
Then $v(m)$ is $(s)$-bounded too.
\end{proposition} 
{\bf Proof:} Let $(a_{t,l})_{t,l}$ be a regulator related with
$(s)$-boundedness of $m$. If the thesis of the
proposition is not true, then we find a disjoint sequence $(A_h)_h$
in ${\mathcal L}$ and a $\varphi\in
{\mathbb{N}}^{\mathbb{N}}$ with 
\begin{eqnarray}\label{hhh}
v(m) (A_h) \not \leq \bigvee_{t=1}^{\infty} a_{t,\varphi(t)}
\quad \text{  for  every  }  h\in\mathbb{N}.
\end{eqnarray}
By (\ref{hhh}) and the properties of the semivariation we find a
sequence $(B_h)_h$ in ${\mathcal L}$, with $B_h \subset A_h$ for 
every $h$ and 
\begin{eqnarray}\label{hhhh}
m(B_h) \not \leq \bigvee_{t=1}^{\infty} a_{t,\varphi(t)}
\quad \text{  for  all  }  h\in\mathbb{N}.
\end{eqnarray}
Since $B_h \subset A_h$ for every $h$, then the sequence $(B_h)_h$
is disjoint, and so (\ref{hhhh}) contradicts $(s)$-boundedness of 
$m$. This ends the proof. $\quad \Box$

Analogously as in Proposition \ref{3434} it is possible to demonstrate 
the following
\begin{proposition}\label{3435}
Let $m_j:{\mathcal L}\to R^+$, $j\in\mathbb{N}$, be a sequence of 
equibounded and uniformly $(s)$-bounded set functions.
Then the set functions $v(m_j):{\mathcal L} \to R^+$, 
$j\in\mathbb{N}$, are
equibounded and uniformly $(s)$-bounded too.
\end{proposition}
The next result extends \cite[Lemma 2 (b)]{SAEKI} to the lattice group setting.
\begin{proposition}\label{3535}
Let $m:{\mathcal L} \to R^+$ be a $k$-triangular set function,
continuous from above at $\emptyset$.
Then $v(m)$ is continuous from above at $\emptyset$.
\end{proposition}  
{\bf Proof:} Let $(a_{t,l})_{t,l}$ be a regulator, associated 
with continuity from above at $\emptyset$ of $m$. To
get the assertion, it will be enough to
prove that $v(m)$ is continuous from above at $\emptyset$
with respect to the $(D)$-sequence $((k+1)a_{t,l})_{t,l}$.
If not, then there are a decreasing
sequence $(A_n)_n$ in ${\mathcal L}$ and an element $\varphi \in
{\mathbb{N}}^{\mathbb{N}}$ with 
$\displaystyle{\bigcap_{n=1}^{\infty} A_n=\emptyset}$ and
\begin{eqnarray}\label{varphisemiv}
v(m) (A_n) \not \leq (k+1) \bigvee_{t=1}^{\infty}  a_{t,\varphi(t)}.
\end{eqnarray}
Now, proceeding by induction, similarly as in 
\cite[Lemma 2 (b)]{SAEKI}, let $D_0 = 0$, $q_0=1$, and suppose that
$q_1 < \ldots < q_{n-1} \in \mathbb{N}$ and the pairwise disjoint sets 
$D_1, \ldots, D_{n-1} \in {\mathcal L}$ have
been given. By (\ref{varphisemiv}) and the 
properties of the semivariation we find a set $B_n \subset 
A_{q_{n-1}}$ with
\begin{eqnarray}\label{varphisemiv2}
m (B_n) \not \leq (k+1)\bigvee_{t=1}^{\infty} a_{t,\varphi(t)}.
\end{eqnarray} By continuity from above at $\emptyset$ of $m$,
for every $\psi\in {\mathbb{N}}^{\mathbb{N}}$ there is 
$\overline{q}\in\mathbb{N}$ (depending on $\psi$ and $n$) with
\begin{eqnarray*}
m (B_n\cap A_q) \leq \bigvee_{t=1}^{\infty} a_{t,\psi(t)}
\quad \text{  whenever  } q \geq \overline{q},
\end{eqnarray*} and thus there exists $q_n> q_{n-1}$ with
\begin{eqnarray}\label{varphisemiv0}
m (B_n\cap A_{q_n}) \leq \bigvee_{t=1}^{\infty} a_{t,\psi(t)}.
\end{eqnarray} Let now $D_n:=B_n \setminus A_{q_n}$. We get 
\begin{eqnarray}\label{varphisemiv1}
m (D_n) \not \leq (k+1) \bigvee_{t=1}^{\infty} a_{t,\varphi(t)}.
\end{eqnarray}
Indeed, if (\ref{varphisemiv1}) is not true, then from 
(\ref{varphisemiv0}), $k$-triangularity of $m$ and
weak $\sigma$-distributivity of $R$ we get
\begin{eqnarray}\label{varphisemiv3}
m (B_n) &\leq& m(D_n) + k \, m(B_n \cap A_{q_n}) \leq
(k+1) \bigvee_{t=1}^{\infty} a_{t,\varphi(t)} + k 
\bigvee_{t=1}^{\infty} a_{t,\psi(t)} \leq \\ &\leq&
(k+1) \bigvee_{t=1}^{\infty} a_{t,\varphi(t)} + k 
\bigwedge_{\psi\in {\mathbb{N}}^{\mathbb{N}}}
\Bigl(\bigvee_{t=1}^{\infty} a_{t,\psi(t)} \Bigr) =
(k+1) \bigvee_{t=1}^{\infty} a_{t,\varphi(t)}, \nonumber
\end{eqnarray}
which contradicts (\ref{varphisemiv2}). Since $D_n
\subset A_{q_n-1} \setminus A_{q_n}$, then 
we obtain that the $D_n$'s are pairwise disjoint
and satisfy (\ref{varphisemiv1}). As ${\mathcal L}$ satisfies
property $(E)$, possibly taking a suitable subsequence
of $(D_n)_n$, without loss of generality we can assume that 
${\mathcal L}$ contains the $\sigma$-algebra generated by the
$D_n$'s. Since $m$ is $k$-triangular and
continuous from above at $\emptyset$
with respect to the regulator $(a_{t,l})_{t,l}$,  then for every 
$\varphi\in{\mathbb{N}}^{\mathbb{N}}$
there is $\overline{n} \in\mathbb{N}$ with 
\begin{eqnarray*}
m(D_n) \leq m \Bigl( \bigcup_{s=n}^{\infty} D_s \Bigr) + k \, 
m \Bigr( \bigcup_{s=n+1}^{\infty} D_s \Bigl) \leq (k+1)
\bigvee_{t=1}^{\infty} a_{t,\varphi(t)} 
\end{eqnarray*}
whenever $n \geq \overline{n}$, getting a contradiction with 
(\ref{varphisemiv1}). This ends the proof. $\quad \Box$.
\vspace{3mm}

Analogously as in Proposition \ref{3535}, it is possible to
prove the following
\begin{proposition}\label{3535bis}
Let $m_j:{\mathcal L} \to R$, $j\in\mathbb{N}$, be a sequence of
equibounded $k$-triangular 
set functions, uniformly continuous from above at $\emptyset$.
Then the set functions $\mu_j:=v(m_j)$, $j\in\mathbb{N}$, are 
uniformly continuous from above at $\emptyset$.
\end{proposition}

We now give an example of $1$-triangular and not monotone
set function. 
\begin{example}\label{measuroid} \rm
Let $R=\mathbb{R}$, and $\mu:{\mathcal P}(\mathbb{N})\to \mathbb{R}$ be
defined as 
\begin{eqnarray}\label{generalizedharmonic}
\mu(A):=\sum_{n\in A} \frac{(-1)^n}{n^2}, \qquad A \subset 
\mathbb{N}.
\end{eqnarray} 
Since the series $\displaystyle{\sum_{n=1}^{\infty}
\frac{(-1)^n}{n^2}}$ is absolutely
convergent, then $\mu$ is a
countably additive measure (see also \cite[Proposition
2.15]{interchange}). 
From this and \cite[Example (a) and 
Lemma 2 (b)]{SAEKI} it follows that 
the set function $m:{\mathcal 
P}(\mathbb{N})\to \mathbb{R}$, 
defined as $m(A):=|\mu(A)|$, is continuous from above
at $\emptyset$. Moreover, by 
construction, $m$ is $1$-triangular.  However, $m$ is not 
monotone: indeed, $\displaystyle{
m(\{1,3\})=\frac{10}{9} >\frac{31}{36}=m(\{1,2,3\})}$. 
\end{example}
\begin{definition}  \rm 
A topology $\tau$ on 
${\mathcal L}$ is a \textit{Fr\'{e}chet-Nikod\'{y}m topology} iff
the functions $(A,B) \mapsto A \Delta B$
and $(A,B) \mapsto A \cap B$ from ${\mathcal L} \times
{\mathcal L}$ (endowed with
the product topology) to ${\mathcal L}$ are continuous,
and for any $\tau$-neighborhood $V$ of
$\emptyset$ in ${\mathcal L}$ there exists a
$\tau$-neighborhood $U$ of $\emptyset$ in
${\mathcal L}$ such that, if
$E\in {\mathcal L}$ is contained in some suitable
element of $U$, then $E \in V$
(see also \cite{DRE2}).
\end{definition}
\begin{definitions}\label{regular} \rm
(a) Let $\tau$ be a Fr\'{e}chet-Nikod\'{y}m topology.
A set function $m:{\mathcal L}\to R^+$ is said to
be \textit{$\tau$-continuous} on ${\mathcal L}$ iff 
it is $(s)$-bounded on ${\mathcal L}$ and for each decreasing 
sequence $(H_n)_n$ in ${\mathcal L}$ with 
$\tau$-$\displaystyle{
\lim_n H_n=\emptyset}$ 
we get $\displaystyle{
(D)\lim_n m(H_n)=0}$ with respect to a single regulator
(see also \cite{DRE2}). 

(b) The set functions $m_j:{\mathcal L}\to
R^+$, $j\in \mathbb{N}$, are \textit{uniformly
$\tau$-continuous} on ${\mathcal L}$ iff they are uniformly 
$(s)$-bounded and for every decreasing sequence 
$(H_n)_n$ in ${\mathcal L}$ with 
$\tau$-$\displaystyle{\lim_n H_n=\emptyset}$ 
we get
$\displaystyle{(D)\lim_n \Bigl(\bigvee_j\, m_j
(H_n)\Bigr)=0}$ with respect to a single regulator.


\end{definitions}
\section{The main results}

We begin with the following proposition, which will be useful
in the sequel.
\begin{proposition}\label{propertyu}
Let $R$ be any Dedekind complete
lattice group, $x\in R$, $(x_n)_n$ be
any sequence in $R$, such that
\begin{description}
\item [{\rm \ref{propertyu}.1)}]
every subsequence $(x_{n_q})_q$ of $(x_n)_n$
has a sub-subsequence $(x_{n_{q_r}})_r$,
$(D)$-convergent to $x$ with respect to a
single $(D)$-sequence $(a_{t,l})_{t,l}$.

Then $\displaystyle{(D)\lim_n x_n =x}$ 
with respect to $(a_{t,l})_{t,l}$.
\end{description}
\end{proposition}
{\bf Proof:} If we deny the thesis, then there exist $\varphi \in
{\mathbb{N}}^{\mathbb{N}}$ and a strictly increasing sequence
$(n_q)_q$ with $\displaystyle{|x_{n_q}-x|\not\leq
\bigvee_{t=1}^{\infty} a_{t,\varphi(t)}}$ for each
$q\in\mathbb{N}$. So any subsequence of $(x_{n_q})_q$ does 
not $(D)$-converge to $x$ with respect to $(a_{t,l})_{t,l}$,
getting a contradiction with \ref{propertyu}.1).
This ends the proof. $\quad \Box$
\begin{remarks} \rm (a) With similar techniques as in Proposition
\ref{propertyu}, it is possible to prove an analogous result also 
when it is dealt with $(O)$-convergence.

(b) Observe that, in general, Theorem \ref{propertyu} 
is not true without convergence of sub-sub\-se\-quen\-ces 
with respect to a single regulator. Indeed, if $\nu$ is the Lebesgue 
measure on $[0,1]$,
${\mathcal M}$ is the $\sigma$-algebra of all $\nu$-measurable
subsets of $[0,1]$ and $R=L^0([0,1],{\mathcal M},\nu)$ is the space 
of all $\nu$-measurable real-valued functions with identification 
up to $\nu$-null sets, then both order- and $(D)$-convergence 
coincide with convergence almost everywhere. Let $(x_n)_n$ 
be a sequence in $R$, convergent in measure but 
not almost everywhere to $0$: note that each subsequence 
$(x_{n_q})_q$ of $(x_n)_n$
has a sub-subsequence $(x_{n_{q_r}})_r$,
convergent almost everywhere to $0$ (see also \cite{bdLAP}). 
\end{remarks}
Now, using the sliding hump technique,
we prove the following Brooks-Jewett-type theorem, which extends 
\cite[Theorem 5.4]{myvhs},
\cite[Theorem 2.6]{CANDELOROPALERMO},
\cite[Theorem 1]{papvhs} and \cite[Theorem 11.10]{PAPBOOK}. 
\begin{theorem}\label{bjktriang} (BJ)
Let 
${\mathcal E}\subset {\mathcal L}$ be a lattice,
satisfying property $(E)$, and $m_j:{\mathcal L}
\to R$, $j\in\mathbb{N}$, be a
sequence of equibounded
$k$-triangular set functions, 
whose restrictions on ${\mathcal E}$
are $(s)$-bounded on ${\mathcal E}$. If the
limit $\displaystyle{m_0(E):=(D)\lim_j m_j(E)}$ exists in
$R$ for every $E\in{\mathcal E}$ with respect to a single
regulator, then the $m_j$'s are uniformly $(s)$-bounded on
${\mathcal E}$, and $m_0$ is $k$-triangular and $(s)$-bounded.
\end{theorem}
{\bf Proof:} Let $\displaystyle{
u:=\bigvee_{j
\in \mathbb{N}}v(m_j)(G)}$. Note that $u\in R$, thanks to
equiboundedness of the $m_j$'s. 
For each $j \in \mathbb{N}$ let
$(a^{(j)}_{t,l})_{t,l}$ be a $(D)$-sequence related with
$(s)$-boundedness of $m_j$. By the Fremlin lemma \ref{fremlin}
there exist a $(D)$-sequence $(a_{t,l})_{t,l}$ with
\begin{eqnarray*}
\displaystyle{u \wedge \Bigl( \bigvee_q
\Bigl(\sum_{j=1}^{q}
\Bigl( \bigvee_{t=1}^{\infty} a^{(j)}_{t,
\varphi(t+j)}\Bigr)
\Bigr) \Bigr)\leq \bigvee_{t=1}^{\infty} a_{t,
\varphi(t)}}
\end{eqnarray*} and a regulator $(\alpha_{t,l})_{t,l}$ with
\begin{eqnarray}\label{argue}
\displaystyle{u \wedge \Bigl( \bigvee_q
\Bigl(\sum_{j=1}^{q}
\Bigl( \bigvee_{t=1}^{\infty} a_{t,
\varphi(t+j+1)}\Bigr)
\Bigr) \Bigr)\leq \bigvee_{t=1}^{\infty} \alpha_{t,
\varphi(t)}}
\end{eqnarray}
for every $\varphi \in
\mathbb{N}^{\mathbb{N}}$.
Let $(b_{t,l})_{t,l}$ be a regulator associated with pointwise
$(D)$-convergence of the $m_j$'s on ${\mathcal E}$. 
We now prove that the regulator $(c_{t,l})_{t,l}$,
defined by $c_{t,l}=2(k+2)^2 (a_{t,l}+b_{t,l})$,
$t$, $l\in\mathbb{N}$, satifies the condition of uniform 
$(s)$-boundedness of the $m_j$'s.
Otherwise we find a disjoint sequence $(C_n)_n$ in ${\mathcal E}
$ and a function $\varphi\in{\mathbb{N}}^{\mathbb{N}}$ such that 
for each $n$ there are $q_n\geq n$ and $j_n\in\mathbb{N}$ with
\begin{eqnarray}\label{absurdhypotheses}
m_{j_n}(C_{q_n}) \not \leq \bigvee_{t=1}^{\infty} c_{t,\varphi(t)}.
\end{eqnarray}
We claim that both the $q_n$'s and the $j_n$'s 
can be taken strictly increasing. Indeed, if
we have found $q_1 < \ldots < q_{n-1}$ and  
$j_1 < \ldots < j_{n-1}$ satisfying (\ref{absurdhypotheses}), then
at the $n$-th step,
by virtue of $(s)$-boundedness of the $m_j$'s,
there is an integer $i_n > q_{n-1}$ with 
\begin{eqnarray}\label{j0j}
{\displaystyle
m_j(C_i) \leq \bigvee_{t=1}^{\infty} c_{t,\varphi(t)}}
\quad   \text{  for  each  } j=1, \ldots, j_{n-1} \text{  and  }
i \geq i_n.
\end{eqnarray} 
By (\ref{absurdhypotheses}), in correspondence with $i_n$ there are
$q_n > i_n$ and $j_n\in\mathbb{N}$, fulfilling 
(\ref{absurdhypotheses}). Since $q_n > i_n$, from (\ref{j0j})
it follows that $j_n > j_{n-1}$, getting the claim.

For each $n\in\mathbb{N}$, let $H_n:=C_{q_n}$. By property $(E)$,
passing to subsequences, without loss of generality we can 
assume that $\displaystyle{\bigcup_{n\in P} H_n
\in {\mathcal E}}$ for every $P\subset \mathbb{N}$.
By $(s)$-boun\-ded\-ness of $m_{j_1}$, we find
two sets $P_1\subset \mathbb{N}$ and $F_1
\in{\mathcal E}$, with $\displaystyle{
F_1=\bigcup_{l\in P_1} H_l}$ 
and $\displaystyle{
v(m_{j_1})(F_1)\leq \bigvee_{t=1}^{\infty} a_{t,\varphi(t+1)}}$.
Put $l_1=\min P_1$: in correspondence with
$H_{l_1}$ there is $n_1>l_1$ with
$\displaystyle{
|m_i(H_{l_1})-m_j(H_{l_1})|\leq \bigvee_{t=1}^{\infty}2\, b_{t,
\varphi(t)}}$ for each $i$, $j\geq n_1$.

By $(s)$-boun\-ded\-ness of $m_{j_{l_1}}$, there are
an infinite set $P_2\subset P_1$, with $l_2=\min
P_2 >l_1$, and a set $F_2 \in{\mathcal E}$, with $\displaystyle{
F_2=\bigcup_{l\in P_2} H_l}$ and $\displaystyle{
v(m_{j_{l_1}})(F_2) \leq 
\bigvee_{t=1}^{\infty} a_{t,\varphi(t+2)}}$. In
correspondence with
$H_{l_2}$ there exists a natural number
$n_2>l_2$ with $\displaystyle{
|m_i(H_{l_2})-m_j(H_{l_2})|\leq \bigvee_{t=1}^{\infty} 2 \,b_{t,\varphi(t)}}$ for any
$i$, $j\geq n_2$. Proceeding by induction, we find
two decreasing sequences $(P_h)_h$, $(F_h)_h$ of infinite
sets and two strictly
increasing sequences $(n_h)_h$, $(l_h)_h$ in
$\mathbb{N}$, satisfying the following
conditions, for every $h\in\mathbb{N}$:
\begin{itemize}
\item[\rm {\ref{bjktriang}.1.)}] $l_h=\min P_h,\, n_h > l_h, \,
l_{h+1}> n_h+l_h$;
\item[\rm {\ref{bjktriang}.2.)}] $\displaystyle{
F_h=\bigcup_{l\in P_h} H_l}$;
\item[\rm {\ref{bjktriang}.3.)}] $\displaystyle{
v(m_{j_{l_h}})(F_{h+1}) \leq \bigvee_{t=1}^{\infty} a_{t,
\varphi(t+h+1)} \leq 
\bigvee_{t=1}^{\infty} \alpha_{t,\varphi(t)}}$;
\item[\rm {\ref{bjktriang}.4.)}] $\displaystyle{
|m_i(H_{j_{l_h}})-m_j(H_{j_{l_h}})|\leq \bigvee_{t=1}^{\infty} 2\,
b_{t,\varphi(t)}}$ for each $i$, $j\geq n_h$.
\end{itemize}
Put $m^{\prime}_h=m_{j_{l_h}}$ and
$H^{\prime}_h=H_{l_h}$, $h\in\mathbb{N}$.
Let $h_r:=2\,r$, $r \in \mathbb{N}$, and
$F:=\displaystyle{\bigcup_{r=1}^{\infty} H^{\prime}_{h_r}}$. Note
that $\displaystyle{\bigcup_{s=r+1}^{\infty} H^{\prime}_{h_s}
 \subset F_{h_r +1}}$ 
and $l_{h_r-1}>n_{h_r-2}\geq n_{h_{r-1}}$
for each $r\in\mathbb{N}$ (see also \cite{CANDELOROPALERMO}).
So we have $$|m^{\prime}_{h_r}
(H^{\prime}_{h_s})-m^{\prime}_{h_r-1}
(H^{\prime}_{h_s})|\leq \bigvee_{t=1}^{\infty} 2 \, 
b_{t,\varphi(t)} \quad \text{for   every  }
s=1, \ldots, r-1.$$ By pointwise convergence of the $m_j$'s with
respect to the regulator $(b_{t,l})_{t,l}$, we find a positive integer
$\overline{r}$ with $$\displaystyle{
m^{\prime}_{h_r}(F) \leq \bigvee_{t=1}^{\infty} b_{t,\varphi(t)}}
\quad \text{  whenever  } r\geq \overline{r}.$$ Since 
$F_{h_r}\subset
F_{h_r-1}$, we have $$m^{\prime}_{h_r}\Bigl(
\bigcup_{s=r+1}^{\infty} 
H^{\prime}_{h_s}\Bigr)\leq\bigvee_{t=1}^{\infty}
\alpha_{t,\varphi(t)} \text{  and  }
m^{\prime}_{h_r-1}\Bigl(
\bigcup_{s=r+1}^{\infty} 
H^{\prime}_{h_s}\Bigr)\leq \bigvee_{t=1}^{\infty} 
\alpha_{t,\varphi(t)}.$$  
From this and Proposition \ref{kfinitetriangularity} we get
\begin{eqnarray*}
m^{\prime}_{h_r}
\Bigl( \bigcup_{s=1}^{r-1} H^{\prime}_{h_s}\Bigr)&\leq&
m^{\prime}_{h_r}
(H^{\prime}_{h_1})+ k \sum_{s=2}^{r-1} m^{\prime}_{h_r}
(H^{\prime}_{h_s})\leq \\&\leq&
k \Bigl( u \wedge \Bigl(\bigvee_q
\Bigl(\sum_{h=1}^q \Bigl(\bigvee_{t=1}^{\infty} a_{t,
\varphi(t+h)}
\Bigr) \Bigr) \Bigr) \Bigr)\leq k\bigvee_{t=1}^{\infty}
\alpha_{t,\varphi(t)},\\
m^{\prime}_{h_r-1}
\Bigl( \bigcup_{s=1}^{r-1} H^{\prime}_{h_s}\Bigr)&\leq&
m^{\prime}_{h_r-1}
(H^{\prime}_{h_1})+ k \sum_{s=2}^{r-1} 
m^{\prime}_{h_r-1}
(H^{\prime}_{h_s})\leq \\ &\leq&
k \Bigl( u \wedge \Bigl(\bigvee_q
\Bigl(\sum_{h=1}^q \Bigl(\bigvee_{t=1}^{\infty} a_{t,
\varphi(t+h)}
\Bigr) \Bigr) \Bigr) \Bigr)\leq k\bigvee_{t=1}^{\infty}
\alpha_{t,\varphi(t)}. 
\end{eqnarray*}
From this and (\ref{kntriangularitybis}) 
used with $q=3$, $E_1=H^{\prime}_{h_r}$, 
$\displaystyle{E_2=\bigcup_{s=r+1}^{\infty} 
H^{\prime}_{h_s}}$, $\displaystyle{E_3= \bigcup_{s=1}^{r-1} 
H^{\prime}_{h_s}}$, since
$F=E_1\cup E_2 \cup E_3$ and  the $m_j$'s are
$k$-triangular, we obtain 
\begin{eqnarray*}
m^{\prime}_{h_r}(H^{\prime}_{h_r})&\leq& k(k+2)
\Bigl(\bigvee_{t=1}^{\infty} (\alpha_{t,\varphi(t)}+b_{t,\varphi(t)})
\Bigr),
\\m^{\prime}_{h_r-1}(H^{\prime}_{h_r})&\leq&  k(k+2)
\Bigl( \bigvee_{t=1}^{\infty} (\alpha_{t,\varphi(t)}+b_{t,\varphi(t)})
\Bigr), \end{eqnarray*} and hence
\begin{eqnarray}\label{19} |m^{\prime}_{h_r}(H^{\prime}_{h_r}) -
m^{\prime}_{h_r-1}(H^{\prime}_{h_r})|\leq2 \, k(k+2)
\Bigl(\bigvee_{t=1}^{\infty}(\alpha_{t,\varphi(t)}+b_{t,\varphi(t)})
\Bigr).\end{eqnarray} Furthermore, thanks to \ref{bjktriang}.3.), we 
get \begin{eqnarray}\label{20}
m^{\prime}_{h_r-1}(H^{\prime}_{h_r}) \leq \bigvee_{t=1}^{\infty} 
\alpha_{t,\varphi(t)}. \end{eqnarray} From (\ref{19}) and (\ref{20})
we obtain $$m^{\prime}_{h_r}(H^{\prime}_{h_r})\leq
2(k+2)^2 \Bigl(\bigvee_{t=1}^{\infty}
(a_{t,\varphi(t)}+b_{t,\varphi(t)})\Bigr),$$ which
contradicts (\ref{absurdhypotheses}). Thus the $m_j$'s are uniformly
$(s)$-bounded on ${\mathcal E}$. From this it is not difficult to
deduce that $m_0$ is $k$-triangular and $(s)$-bounded on 
${\mathcal E}$. $\quad\Box$ \vspace{3mm}

The following result will be useful to prove
our versions of the Vitali-Hahn-Saks and Nikod\'{y}m theorem,
and extends \cite[Corollaries 3.5 and 5.5]{myvhs} and
\cite[Lemma 3.13]{CANDELOROLETTA}.
\begin{lemma}\label{crucialpointgl}
Let 
${\mathcal G}$ and ${\mathcal H}$ be two sublattices of 
${\mathcal L}$, such that the complement of every 
element of ${\mathcal H}$ belongs to ${\mathcal G}$,
$m_j:{\mathcal L} \to R$, $j \in
\mathbb{N}$, be a sequence of $k$-triangular set functions,
uniformly $(s)$-bounded on
${\mathcal G}$. Fix $W \in {\mathcal H}$ and a decreasing
sequence $(H_n)_n$ in ${\mathcal G}$, with
$W\subset H_n$ for each $n\in \mathbb{N}$. If
\begin{eqnarray}\label{hypothesisgl}
(D)\lim_n \Bigl( \bigvee_{A\in {\mathcal G}, A \subset 
H_n \setminus W}m_j(A)\Bigr)=\bigwedge_n \Bigl( \bigvee_{A\in 
{\mathcal G}, A \subset H_n \setminus W} m_j(A)\Bigr)=0  
\text{  for every  } j \in \mathbb{N}
\end{eqnarray}
with respect to a single
$(D)$-sequence $(a_{t,l})_{t,l}$, then
\begin{eqnarray*}
(D)\lim_n \Bigl(\bigvee_j \Bigl( \bigvee_{A\in {\mathcal G}, A \subset 
H_n \setminus W}m_j(A)\Bigr) \Bigr)=\bigwedge_n 
\Bigl(\bigvee_j \Bigl(\bigvee_{A\in {\mathcal G}, A \subset 
H_n \setminus W}m_j(A)\Bigr) \Bigr)=0
\end{eqnarray*}
with respect to $(a_{t,l})_{t,l}$.
\end{lemma}
{\bf Proof:} Put ${\mathcal W}:=\{A\in{\mathcal G}:A
\cap W=\emptyset\}$. By $k$-triangularity of $m_j$, for
any $A \in {\mathcal
W}$,  $j$, $q\in\mathbb{N}$, we have
\begin{eqnarray*}
m_j(A\setminus H_q) &\leq& k \, m_j(A \cap H_q) + m_j(A),
\\ m_j(A)&\leq& m_j(A \setminus H_q) + k \, m_j(A \cap H_q),
\end{eqnarray*} and hence
\begin{eqnarray}\label{ktriang1gl}
0 &\leq& |m_j(A)- m_j(A \setminus H_q)| =
\\&=&(m_j(A)-m_j(A \setminus H_q) ) \vee 
(m_j(A \setminus H_q) - m_j(A)) \nonumber\leq k \, m_j(A \cap H_q).
\end{eqnarray}
From (\ref{hypothesisgl}) and (\ref{ktriang1gl}), since $A \cap H_q
\subset H_{q-1} \setminus W$ for every $q\in \mathbb{N}$, we get
\begin{eqnarray}\label{qugl}
m_j(A)=(D)\lim_q \, m_j(A\setminus H_q)\quad
\text{  for each  } A \in {\mathcal W} \text{  and  }
j\in\mathbb{N}.
\end{eqnarray}
If we deny the thesis of the lemma, then there
is a $\varphi\in\mathbb{N}$ such that for every $r \in
\mathbb{N}$ there are $j$, $n \in \mathbb{N}$
with $n>r$ and $A\in{\mathcal G}$ with $A
\subset H_n\setminus W$, $\displaystyle{m_j(A)\not\leq 
\bigvee_{t=1}^{\infty} a_{t,\varphi(t)}}
$, and thus, thanks to (\ref{qugl}),
\begin{eqnarray*}\label{absurdhypothesisgl}
m_j(A\setminus H_q)\not\leq\bigvee_{t=1}^{\infty} a_{t,\varphi(t)}
\end{eqnarray*}
for $q$ large enough.

At the first step, we find a set
$A_1\in{\mathcal G}$ and three integers $n_1 >1$,
$j_1\in\mathbb{N}$ and $q_1 >
\max\{j_1,n_1\}$, with $A_1 \subset H_{n_1}
\setminus W$, $\displaystyle{
m_{j_1}(A_1)\not\leq\bigvee_{t=1}^{\infty} a_{t,\varphi(t)}}$
and $\displaystyle{
m_{j_1}(A_1\setminus H_{q_1})\not\leq\bigvee_{t=1}^{\infty} a_{t,\varphi(t)}}$.
From (\ref{hypothesisgl}), in correspondence with
$j=1, 2, \ldots, j_1$ there exists $h_1>q_1$ with
\begin{eqnarray}\label{ngl}
m_j(A)\leq \bigvee_{t=1}^{\infty} a_{t,\varphi(t)}
\end{eqnarray} whenever $n\geq h_1$ and
$A\subset H_n\setminus W$.

At the second step, there are $A_2\in{\mathcal G}$,
$n_2>h_1$, $j_2\in\mathbb{N}$ and $q_2 >
\max\{j_2,n_2\}$, with $A_2 \subset H_{k_2}
\setminus W$ and
\begin{eqnarray}\label{q2absurdgl}
m_{j_2}(A_2)\not\leq\bigvee_{t=1}^{\infty} a_{t,\varphi(t)};\quad m_{j_2}
(A_2 \setminus H_{q_2}) \not \leq \bigvee_{t=1}^{\infty} a_{t,\varphi(t)}.
\end{eqnarray}
From (\ref{ngl}) and (\ref{q2absurdgl}) it follows
that $j_2>j_1$.

Proceeding by induction, we find a sequence $(A_r)_r$ in
${\mathcal G}$ and three strictly increasing
sequences in $\mathbb{N}$, $(h_r)_r$, $(j_r)_r$,
$(q_r)_r$, with $q_r >n_r >q_{r-1}$ for any $r\geq
2$; $q_r>j_r$, $A_r\subset H_{n_r}\setminus W$,
$\displaystyle{m_{j_r}(A_r\setminus H_{q_r})\not \leq
\bigvee_{t=1}^{\infty} a_{t,\varphi(t)}}$ for each $r\in\mathbb{N}$. 
But this is impossible, since the sets $A_r\setminus
H_{q_r}$, $r\in\mathbb{N}$, are disjoint elements of ${\mathcal G}$
and the measures $m_j$, $j\in\mathbb{N}$, are
globally uniformly $(s)$-bounded on ${\mathcal G}$ with
respect to $(a_{t,l})_{t,l}$. 
This ends the proof. $\quad \Box$ \vspace{3mm}

As a consequence of Theorem \ref{bjktriang} and
Lemma \ref{crucialpointgl}, we deduce
the following Nikod\'{y}m-type theorem.  
\begin{theorem}\label{nikodym}  (N) Let
${\mathcal L}$ 
satisfy property $(E)$,
$m_j:{\mathcal L} \to R$,
$j\in\mathbb{N}$, be a sequence of equibounded 
$k$-triangular set functions, continuous from above at
$\emptyset$. Suppose that the limit 
$\displaystyle{m_0(E):=(D)\lim_j m_j(E)}$
exists in $R$ for each $E\in {\mathcal L}$
with respect to a single regulator. 

Then the $m_j$'s are uniformly continuous from above at 
$\emptyset$, and $m_0$ is $k$-triangular and continuous 
from above at $\emptyset$ on ${\mathcal L}$.
\end{theorem}
{\bf Proof:} First of all, observe that the $m_j$'s are 
$(s)$-bounded on ${\mathcal L}$: indeed, 
if $(C_h)_h$ is any disjoint sequence in ${\mathcal L}$, then we have
\begin{eqnarray}\label{-22}
0 \leq v(m_j)(C_h) \leq v(m_j) \Bigl( \bigcup_{n=h}^{\infty} C_n \Bigr)
\quad \text{  for every  } j, h \in\mathbb{N}.
\end{eqnarray}
From (\ref{-22}), continuity from above at $\emptyset$
of $m_j$ and monotonicity of $v(m_j)$ we get
$\displaystyle{(D)\lim_h v(m_j) (C_h)=0}$, namely $(s)$-boundedness
of $m_j$. From this and Theorem \ref{bjktriang} used with ${\mathcal 
E}={\mathcal L}$ we deduce uniform $(s)$-boundedness of the 
$m_j$'s.

Choose arbitrarily $j\in\mathbb{N}$ and
a decreasing sequence $(H_n)_n$ in
${\mathcal L}$ with $\displaystyle{\bigcap_{n=1}^{\infty} 
H_n=\emptyset}$. Since $m_j$ is continuous from above at 
$\emptyset$, we 
have $$\displaystyle{(D)\lim_n v(m_j)(H_n)= \bigwedge_n 
v(m_j) (H_n)=0}.$$
By Lemma \ref{crucialpointgl} used with ${\mathcal G}={\mathcal
H}={\mathcal L}$ and $W=\emptyset$ and taking into account uniform 
$(s)$-boundedness of the $m_j$'s, we have 
$$\displaystyle{(D)\lim_n \Bigl( \bigvee_j v(m_j)(H_n)\Bigr)
=\bigwedge_n \Bigl( \bigvee_j 
v(m_j)(H_n)\Bigr)=0}.$$ By arbitrariness of the sequence
$(H_n)_n$, we get
uniform continuity from above at $\emptyset$ of the $m_j$'s.
From this it follows easily that $m_0$ is $k$-triangular and
continuous from above at $\emptyset$.
$\quad \Box$
\vspace{3mm}

Analogously as Theorem \ref{nikodym}, it is possible to prove
the following Vitali-Hahn-Saks-type theorem.
\begin{theorem}\label{vhs} (VHS)
Let $R$, $G$, ${\mathcal L}$ be as in Theorem \ref{nikodym},
$\tau$ be a Fr\'{e}chet-Nikod\'{y}m topology on ${\mathcal L}$,
$m_j:{\mathcal L} \to R$,
$j\in\mathbb{N}$, be a sequence of equibounded
$\tau$-continuous $k$-triangular set functions. Let
$\displaystyle{m_0(E):=(D)\lim_j m_j(E)}$
exist in $R$ for every $E\in {\mathcal L}$
with respect to a single regulator. 

Then the $m_j$'s are uniformly $\tau$-continuous on ${\mathcal L}$,
and $m_0$ is $k$-triangular and $\tau$-continuous on ${\mathcal L}$.
\end{theorem}
We now prove a Schur-type theorem, which extends 
\cite[Corollary 5.6]{myvhs}. 
\begin{theorem}\label{schur} (S) Let $R$ be any
Dedekind complete and weakly $\sigma$-distributive 
lattice group, $m_j:{\mathcal P}(\mathbb{N}) \to R^+$, $j \in
\mathbb{N}$, be a sequence of
equibounded, continuous from above at $\emptyset$ and
$k$-triangular set functions, and let there
exist a set function $m_0:{\cal P}(\mathbb{N})
\to R$ with $\displaystyle{(D)\lim_j \, m_j(E)=m_0(E)}$
for every $E\subset \mathbb{N}$ with respect to a single
$(D)$-sequence. 

Then the $m_j$'s are  
uniformly continuous from above at $\emptyset$ and
$m_0$ is $k$-triangular and 
continuous from above at $\emptyset$. Furthermore we get
\begin{eqnarray}\label{estimate1}
(D)\lim_j \Bigl(\bigvee_{E\subset \mathbb{N}} 
|m_j(E)-m_0(E)|\Bigr)=0.
\end{eqnarray}
\end{theorem}
{\bf Proof:}  Uniform continuity from above at $\emptyset$ of the 
$m_j$'s, $k$-triangularity and continuity from above at $\emptyset$ 
of $m_0$ follow from Theorem \ref{nikodym}.
In particular, there exists a regulator $(a_{t,l})_{t,l}$ such that
for each $\varphi \in {\mathbb{N}}^{\mathbb{N}}$
there is $h_0 \in \mathbb{N}$ with
\begin{eqnarray}\label{k01}
m_j(E\cap [h_0,+\infty[)\leq \bigvee_{t=1}^{\infty} a_{t,\varphi(t)}
\end{eqnarray} for any $j\in \mathbb{N} \cup\{0\}$ 
and $E\subset \mathbb{N}$. 
Let now $(b_{t,l})_{t,l}$ be a $(D)$-sequence associated 
with (pointwise) convergence of $(m_j)_j$ to $m_0$. 
We will prove that the regulator $((2\,k+2)(a_{t,l}+b_{t,l}))_{t,l}$,
$t$, $l\in\mathbb{N}$, satisfies the condition of 
$(D)$-limit in (\ref{estimate1}).
In correspondence with $\varphi$ and $h_0$ as in (\ref{k01})
there exists a positive integer $j_0$ with 
\begin{eqnarray}\label{j01}
|m_j (E\cap [1,h_0-1]) - m_0 (E\cap [1,h_0-1])|\leq  
\bigvee_{t=1}^{\infty} a_{t,\varphi(t)}
\end{eqnarray} whenever $j \geq j_0$ and $E\subset \mathbb{N}$.
Since $m_j$ is $k$-triangular for every $j\in \mathbb{N} \cup \{0\}$,
we get \begin{eqnarray*}
& & -k \, m_j(E\cap [h_0,+\infty[)-k \, m_0(E\cap [h_0,+\infty[) \leq
\\ &\leq&
m_j(E) -m_0(E) -m_j(E\cap [1,h_0-1])+m_0(E\cap [1,h_0-1]) \leq 
\\&\leq&k \,m_j(E\cap [h_0,+\infty[)+k \, m_0(E\cap [h_0,+\infty[),
\end{eqnarray*} 
and hence \begin{eqnarray*}
& & |m_j(E) -m_0(E)|\leq |m_j(E\cap [1,h_0-1])+m_0(E\cap 
[1,h_0-1])| + \\ &+&
k \, m_j(E\cap [h_0,+\infty[) + k \, m_0(E\cap [h_0,+\infty[) \leq 
\bigvee_{t=1}^{\infty} (2\,k+2) (a_{t,\varphi(t)}+b_{t,\varphi(t)})
\end{eqnarray*} for every $j \geq j_0$ and $E\subset \mathbb{N}$,
getting the assertion. $\quad \Box$ 
\normalsize
\vspace{3mm}

We now prove the following Drewnowski-type theorem, which extends
\cite[Lemma 2.3]{CANDELOROLETTA} and 
\cite[Theorem 5.3]{bcdefin} to non-additive lattice group-valued set 
functions.
\begin{theorem} \label{subsigma}
Let $R$ be a super Dedekind complete and weakly 
$\sigma$-distributive lattice group, 
$G$, ${\mathcal L}$ be as in Theorem \ref{nikodym}, and
$m:{\mathcal L}\to R$ be any positive $(s)$-bounded set
function. Then for every disjoint
sequence $(C_n)_n$ in ${\mathcal L}$ there exists a subsequence 
$(C_{n_h})_h$
such that $m$ is continuous from above at $\emptyset$ 
on the $\sigma$-algebra
generated by $(C_{n_h})_h$ in the set  
$\displaystyle{\bigcup_{h=1}^{\infty}C_{n_h}}$.
\end{theorem}
{\bf Proof:} Let $(a_{t,l})_{t,l}$ be a regulator associated with
$(s)$-boundedness of $m$, and set 
$$\Phi:=\Bigr\{\bigvee_{t=1}^{\infty}a_{t,\varphi(t)}: \varphi\in 
\mathbb{N}^{\mathbb{N}}\Bigr\}.$$ Since
$R$ is weakly $\sigma$-distributive, we get $\bigwedge \Phi=0$.
By super Dedekind completeness of $R$ there is a sequence
$(\varphi_l)_l$ in $\mathbb{N}^{\mathbb{N}}$ (which 
without loss of generality we can take 
increasing, similarly as in 
\cite[Theorem 5.3]{bcdefin}), with
$$0= \bigwedge \Phi=\bigwedge \Bigl\{\bigvee_{t=1}^{\infty}a_{t,
\varphi_l(t)}:l\in \mathbb{N}\Bigr\}.$$ Set $\displaystyle{
b_l:=\bigvee_{t=1}^{\infty}a_{t,\varphi_l(t)}}$,
$l\in\mathbb{N}$. It is not difficult to see that $(b_l)_l$ is
an $(O)$-sequence in $R$.
Put now $b_{t,l}:=b_l$,  $t$, $l\in \mathbb{N}$. It is 
readily seen that $(b_{t,l})_{t,l}$ is a regulator
(see also \cite{bonitra}): we will show
that it satisfies the assertion.

Choose arbitrarily any disjoint sequence $(C_n)_n$ in ${\mathcal L}$.
By property $(E)$, there is a subsequence $(H_n)_n$ of $(C_n)_n$,
such that ${\mathcal L}$ contains the $\sigma$-algebra 
generated by the $H_n$'s in the set 
$\displaystyle{\bigcup_{n=1}^{\infty}
H_n}$. Choose any disjoint
sequence $(P_r^1)_r$ of infinite subsets of $\mathbb{N}$,
and for each $r\in\mathbb{N}$ 
define $H_r^1:=\bigcup\{H_n:n\in P_r^1\}$.
The sequence $(H_r^1)_r$ is a disjoint sequence in ${\mathcal L}$, 
and by $(s)$-boundedness of $m$
we find $r_1\in \mathbb{N}$ with
$$v(m)(H_r^1)\leq \bigvee_{i=1}^{\infty}a_{i,\varphi_1(i)}=b_1$$
for all $r\geq r_1$.
We now consider $P^1_{r_1}$, and choose any infinite partition of it
into a sequence of disjoint infinite subsets $(P^2_r)_r$.
For each $r\in\mathbb{N}$ set $\displaystyle{
H^2_r:=\bigcup\{H_n:n\in P^2_r\}}$.
Note that the $H^2_r$'s are pairwise disjoint, and each of them is 
contained in $H^1_{r_1}$ by construction.
Again by $(s)$-boundedness of $m$, there exists an integer 
$r_2>r_1$ with $v(m)(H^2_r)\leq b_2$ for any $r\geq r_2.$ 
Proceeding by induction, we obtain a decreasing sequence 
$(P^l_{r_l})_l$ of
infinite subsets of $\mathbb{N}$, and a corresponding sequence
$(H^l_{r_l})_l$, $\displaystyle{
H^l_{r_l}=\cup\{H_n:n\in P^l_{r_l}\}}$,
satisfying $v(m)(H^l_{r_l})\leq b_l$ for any $l$. Let us denote by
$n_1$ the first element of $P^1_{r_1}$, by $n_2$ the first element
of $P^2_{r_2}$ larger than $n_1,$ and so on. We claim that the
sequence $(H_{n_h})_h$ is the requested one.

Indeed, observe that $H_{n_{h+p}}\subset H^h_{n_h}$ for every $h$,
$p\in \mathbb{N}$. Set $\displaystyle{H^*=\bigcup_{h=1}^{\infty} 
H_{n_h}}$: note that, thanks to property $(E)$, $H^*\in {\mathcal 
L}$. Choose any decreasing sequence $(F_s)_s$ in the 
$\sigma$-algebra
generated in $H^*$ by the sets $H_{n_h}$, with $\displaystyle{
\bigcap_{s=1}^{\infty} F_s=\emptyset}$. For every $s$ there exists 
an integer $j(s)$ with $\displaystyle{F_s\subset \bigcup_{h\geq 
j(s)}H_{n_h}}$. Note that $\displaystyle{\lim_s j(s)=+\infty}$.
Choose now $\varphi\in \mathbb{N}^{\mathbb{N}}$ and pick any 
integer $s_0$ with
$j(s_0)\geq\varphi(1)$. For every $s\geq s_0$ we get $$v(m)(F_s)\leq
v(m)(\bigcup_{h\geq j(s_0)}H_{n_h})\leq
v(m)(H^{j(s_0)}_{r_{j(s_0)}})\leq b_{j(s_0)}\leq b_{\varphi(1)}\leq
\bigvee_{t=1}^{\infty} b_{t,\varphi(t)}. $$
This ends the proof. $\quad \Box$ 
\vspace{3mm}

A consequence of Theorem \ref{subsigma} is the following
\begin{theorem}\label{drewnowski} 
Let $R$, $G$, ${\mathcal L}$ be as in Theorem \ref{subsigma},
$m_j:{\mathcal L} \to R$, $n\in \mathbb{N}$, be
a sequence of equibounded, positive and $(s)$-bounded set 
functions. Then
for every disjoint sequence $(C_h)_h$ in ${\mathcal L}$
there is a subsequence $(C_{n_h})_h$ such that every $m_j$ is 
continuous from above at $\emptyset$ on the $\sigma$-algebra 
generated by $(C_{n_h})_h$ in the set $\displaystyle{
\bigcup_{h=1}^{\infty} C_{n_h}}$.
\end{theorem}
{\bf Proof:} Let $\displaystyle{
u:=\bigvee_{ j
\in \mathbb{N}}v(m_j)(G)}$: by
equiboundedness of the $m_j$'s we get that $u\in R$. 
For every $j \in \mathbb{N}$ let
$(a^{(j)}_{t,l})_{t,l}$ be a $(D)$-sequence related with
$(s)$-boundedness of $m_j$. By the Fremlin lemma \ref{fremlin}
there exist a $(D)$-sequence $(a_{t,l})_{t,l}$ with
\begin{eqnarray*}
\displaystyle{u \wedge \Bigl( \bigvee_q
\Bigl(\sum_{j=1}^{q}
\Bigl( \bigvee_{t=1}^{\infty} a^{(j)}_{t,
\varphi(t+j)}\Bigr)
\Bigr) \Bigr)\leq \bigvee_{t=1}^{\infty} a_{t,
\varphi(t)}}.
\end{eqnarray*} Set now  
$\displaystyle{
\Psi:=\Bigr\{\bigvee_{t=1}^{\infty}a_{t,\varphi(t)}: \varphi\in 
\mathbb{N}^{\mathbb{N}}\Bigr\}}$. Since
$R$ is weakly $\sigma$-distributive, we get $\bigwedge \Psi=0$.
Arguing analogously as in Theorem \ref{subsigma}, we find an
increasing sequence 
$(\varphi_l)_l$ in $\mathbb{N}^{\mathbb{N}}$, with
$$0= \bigwedge \Psi=\bigwedge \Bigl\{\bigvee_{t=1}^{\infty}a_{t,
\varphi_l(t)}:l\in \mathbb{N}\Bigr\}.$$ Set $\displaystyle{
b_l:=\bigvee_{t=1}^{\infty}a_{t,\varphi_l(t)}}$,
$l\in\mathbb{N}$. Then $(b_l)_l$ is
an $(O)$-sequence in $R$, and the double sequence 
$(b_{t,l})_{t,l}$, defined by setting $b_{t,l}:=b_l$,  $t$, $l\in 
\mathbb{N}$, is a $(D)$-sequence. We now prove 
that it fulfils the assertion.

By property $(E)$, the sequence $(C_n)_n$ admits a
subsequence $(H_n)_n$, such that ${\mathcal L}$ contains
the $\sigma$-algebra generated by the $H_n$'s in 
$\displaystyle{\bigcup_{n=1}^{\infty} H_n}$.
By Theorem \ref{subsigma}
there exist an infinite subset $P_1 \subset
\mathbb{N}$ and a positive integer $h_1$ 
with $$\displaystyle{
v(m_1)\Bigr(\bigcup_{j \in
P_1,  j \geq h}H_j\Bigr) \leq b_1 \quad \text{  whenever  }
h \geq h_1.}$$
There are an infinite subset $P_2
\subset P_1$ and a $h_2 > h_1$ with
$$\displaystyle{v(m_2)\Bigl(\bigcup_{j \in
P_2, j \geq h}H_j\Bigr)\leq b_2} \quad \text{  for each  }
h \geq h_2.$$ Without loss of generality, we can and
do suppose $\min P_2 > \min P_1$.

Proceeding by induction, we find a decreasing
sequence $(P_n)_n$ of infinite subsets of
$\mathbb{N}$, a strictly increasing sequence
$(p_n)_n$ in $\mathbb{N}$, with $p_n=\min
P_n$ for all $n$, and a sequence
$(h_n)_n$ with $h_n<
h_{n+1}$ for each $n\in\mathbb{N}$,
and $\displaystyle{v(m_n)\Bigl(\bigcup_{j \in P_n,
j \geq h}H_j\Bigr)\leq b_n}$ for every $h \geq
h_n$. Let $P:=\{p_n: n \in \mathbb{N}\}$,
and set $q_n:=\max \{h_n, p_n\}$. For each
$n$ and $h \geq q_n$, we get
\begin{eqnarray*} 0 \leq 
v(m_n)\Bigl(\bigcup_{j \in P,  j \geq h} H_j\Bigr)
\leq v(m_n)\Bigl(\bigcup_{j \in P_n,  j \geq h} H_j
\Bigr) \leq b_n .
\end{eqnarray*}
From this it follows that $$(D)\lim_n
v(m_n)\Bigl(\bigcup_{j \in P,  j \geq h} H_j\Bigr)=0=\bigwedge_n
v(m_n)\Bigl(\bigcup_{j \in P,  j \geq h} H_j\Bigr)$$
for each $n\in\mathbb{N}$, getting the assertion.
$\quad \Box$ \vspace{3mm}

We now prove the equivalence between Theorems 
$(BJ)$, $(N)$, $(VHS)$ and $(S)$, extending
\cite[Theorem, p. 726]{DRE2}. \vspace{3mm}

$(BJ) \Longrightarrow (VHS)$ See Theorem \ref{vhs}.
\vspace{3mm}

$(VHS) \Longrightarrow (N)$ (see also \cite{bdequiv, DRE2})
Let $\tau$ be the Fr\'{e}chet-Nikod\'{y}m
topology generated by the family of all order continuous 
submeasures, 
that is all $1$-subadditive, increasing and 
continuous from above at $\emptyset$ real-valued set 
functions defined on ${\mathcal L}$
(see also \cite{dobrakov1, 
DRE2}).
If $(H_n)_n$ is any decreasing
sequence in ${\mathcal L}$ with $\tau$-$
\displaystyle{\lim_n H_n=
\emptyset}$ and $\displaystyle{H=
\bigcap_{n=1}^{\infty}H_n}$ $\in {\mathcal L}$,
then $\eta(H)=0$
for every order continuous submeasure
$\eta$, and so $H=\emptyset$. From this and Remark \ref{29} 
it follows that, if $m_j:{\mathcal L} \to
R$, $j \in \mathbb{N}$, is a sequence of set functions,
continuous from above at $\emptyset$, then they are
$\tau$-continuous. By $(VHS)$, they are 
uniformly $\tau$-continuous. From this and Lemma 
\ref{crucialpointgl} it follows that the $m_j$'s are
also uniformly continuous from above at $\emptyset$.
Thus, $(VHS)$ implies $(N)$. \vspace{3mm}

$(N) \Longrightarrow (BJ)$ 
Let $(C_h)_h$ be any disjoint sequence in 
${\mathcal L}$ and $m_j:{\mathcal L} \to R$, $j\in\mathbb{N}$,
be a sequence of equibounded $k$-triangular $(s)$-bounded set
functions. Fix arbitrarily a subsequence $(C_{h_r})_r$ of $(C_h)_h$.
By Theorem \ref{drewnowski} there is a sub-subsequence 
$(C_{h_{r_s}})_s$ of $(C_{h_r})_r$ such that every $m_j$ is
continuous from above at $\emptyset$ on the $\sigma$-algebra
generated by the $C_{h_{r_s}}$'s in the set $\displaystyle{
\bigcup_{s=1}^{\infty} C_{h_{r_s}} }$. From this and $(N)$ it
follows that $\displaystyle{(D)\lim_s \Bigl(\bigvee_j v(m_j)
(C_{h_{r_s}})\Bigr)=0 }$ with respect to a suitable regulator 
$(b_{t,l})_{t,l}$, independent of $(C_{h_r})_r$.
By arbitrariness of $(C_{h_r})_r$ and Proposition 
\ref{propertyu}
we get $$(D)\lim_h \Bigl(\bigvee_j v(m_j)(C_h)\Bigr)=0$$
with respect to $(b_{t,l})_{t,l}$,
that is the assertion. 

$(N) \Longrightarrow (S)$  See Theorem \ref{schur}.

$(S) \Longrightarrow (N)$ 
Let ${\mathcal L}$ and $m_j$, $j\in\mathbb{N}$,
be as in the hypotheses of Theorem \ref{nikodym} and 
$(A_n)_n$ be a 
disjoint sequence in 
${\mathcal L}$.
Choose arbitrarily a subsequence $(C_h)_h$ of $(A_h)_h$. 
By property $(E)$ there is a subsequence $(C_{h_r})_r$ of
$(C_h)_h$, such that ${\mathcal L}$ contains the  
$\sigma$-algebra generated by the $C_{h_r}$'s in 
$\displaystyle{\bigcup_{r=1}^{\infty} C_{h_r}}$. 
For every $j\in\mathbb{N}$ and $A\subset \mathbb{N}$, set 
\begin{eqnarray}\label{mu}
\mu_j(A)=m_j \Bigl( \bigcup_{r\in A} C_{h_r} \Bigr).
\end{eqnarray}
We claim that $\mu_j$ is continuous from above at $\emptyset$
for every $j\in\mathbb{N}$.
Fix any decreasing sequence $(K_n)_n$ in ${\mathcal P}(
\mathbb{N})$.
Without loss of generality, we can suppose that 
$(K_n)_n$ is strictly decreasing, and $K_n \subset [n,+\infty[$
for each $n\in\mathbb{N}$, getting
\begin{eqnarray*}\label{vimes}
0 \leq v(\mu_j) (K_n)\leq v(\mu_j)([n,+\infty[) \leq v(m_j)
(\bigcup_{r=n}^{\infty} C_{h_r})
\end{eqnarray*} for every $j$, $n\in\mathbb{N}$. 
From (\ref{vimes}), Proposition \ref{3535} and continuity from above 
at $\emptyset$ of $m_j$ (with respect to a single regulator 
$(a_{t,l})_{t,l}$, which without loss of generality can be taken 
independent of $(C_h)_h$ and of 
$j$, arguing analogously as in (\ref{argue}) ) we get
$$(D)\lim_n v(m_j)(\bigcup_{r=n}^{\infty} C_{h_r})=0$$
with respect to $(a_{t,l})_{t,l}$, and thus we have also
$\displaystyle{(D)\lim_n v(\mu_j)(K_n)=0}$ 
with respect to $(a_{t,l})_{t,l}$, proving the claim.
Moreover, it is not difficult to check pointwise convergence of the
$\mu_j$'s with respect to a single regulator $(b_{t,l})_{t,l}$. 
By Theorem \ref{schur} the $\mu_j$'s are uniformly 
continuous from above at $\emptyset$.
In particular, we get
\begin{eqnarray*} 0\leq
(D)\lim_n \Bigl( \bigvee_j v(m_j)(C_{h_n}) \Bigr)\leq
(D)\lim_n \Bigl( \bigvee_j \mu_j ([n,+\infty[) \Bigr) =0, 
\end{eqnarray*} namely
\begin{eqnarray}\label{f}
(D)\lim_n \Bigl( \bigvee_j v(m_j)(C_{h_n}) \Bigr)=0
 \end{eqnarray} 
with respect to a regulator, which by construction can be taken
independent of the chosen sequence $(C_h)_h$. By arbitrariness
of $(C_h)_h$, (\ref{f}) and Proposition \ref{propertyu} we get
$$\displaystyle{(D)\lim_h \Bigl( \bigvee_j v(m_j)(A_h) \Bigr)=0}.$$ 
Thus the $m_j$'s are uniformly $(s)$-bounded. From this and 
Lemma \ref{crucialpointgl} we get uniform continuity from above 
at $\emptyset$ of the $m_j$'s, and hence $(N)$.
$\quad \Box$ \vspace{3mm}

{\bf Open problems:} 

(a) Prove some limit theorems for $k$-triangular set functions
with respect to other kinds of $(s)$-boundedness and continuity,
and/or relatively to some other types of convergence.

(b) Investigate some other property of non-additive
set functions and their (semi)variations
(see also \cite{BAN, PAPBOOK}).


\normalsize
\begin{thebibliography}{99}
\bibitem{KADETS} A. Aviles Lopez, B. Cascales
Salinas, V. Kadets and A. Leonov,
The Schur $l_1$ theorem for filters,
\textit{J. Math. Phys.,
Anal., Geom.} \textbf{3} (4) (2007), 383-398.
\bibitem{BAN} A. I. Ban, \textit{Intuitionistic fuzzy sets:
theory and applications}, Nova Science Publ., Inc.,
New York, 2006.
\bibitem{bonitra} A. Boccuto, Egorov property and
weak $\sigma$-distributivity in Riesz spaces,
\textit{Acta Math. (Nitra)} \textbf{6} (2003),
61-66. 
\bibitem{myvhs} A. Boccuto and D. Candeloro,
Uniform $(s)$-boundedness and
con\-ver\-gen\-ce  results for measures with
values in complete
$(\ell)$-groups, \textit{J. Math. Anal. Appl.}
\textbf{265} (2002), 170-194.
\bibitem{bcdefin} A. Boccuto and D. Candeloro,
Convergence and decompositions for
$(\ell)$-group-valued set functions,
\textit{Comment. Math.} \textbf{44} (1) (2004),
11-37.
\bibitem{bdequiv} A. Boccuto and X. Dimitriou,
Ideal limit theorems and their equivalence
in $(\ell)$-group setting, \textit{J. Math.
Research} \textbf{5} (2) (2013), 43-60.
\bibitem{bdsubadditive} A. Boccuto and X. Dimitriou,
Limit theorems for $k$-subadditive lattice group-valued 
capacities in the filter convergence setting, \textit{Tatra Mt. 
Math. Publ.} \textbf{65} (2016), 1-21.
\bibitem{bdbook} A. Boccuto and X. Dimitriou,
\textit{Convergence Theorems for
Lattice Group-Valued Measures}, Bentham Science Publ.,
U. A. E., 2015. 
\bibitem{interchange} A. Boccuto and X. Dimitriou,
Matrix theorems and interchange
for lattice group-valued series in the filter convergence setting,
2015, \textit{Bull. Hellenic Math. Soc.} 
\textbf{59} (2016), 39-55.
\bibitem{bdgreece2016} A. Boccuto and X. Dimitriou,
Limit theorems for lattice group-valued $k$-triangular set 
functions,  \textit{Proceedings of the 33rd  PanHellenic 
Conference on Mathematical Education, Chania, Greece, 4-6 
November 2016} (2016), 1-10.
\bibitem{bdgreece2017} A. Boccuto and X. Dimitriou,
Schur-type Theorems for $k$-Triangular Lattice Group-Valued 
Set Functions with Respect to Filter Convergence, 
\textit{Appl. Math. Sci.} \textbf{11} (57) (2017),  
2825--2833.
\bibitem{bdLAP} A. Boccuto and X. Dimitriou, 
\textit{Non-additive lattice group-valued 
set functions and limit theorems}, Lambert Acad. Publ., 2017.
ISBN 978-613-4-91335-5.
\bibitem{bdpmodena} A. Boccuto, X. Dimitriou
and N. Papanastassiou, Countably additive
restrictions and limit theorems in
$(\ell)$-groups, \textit{Atti Semin. Mat. Fis.
Univ. Modena Reggio Emilia} \textbf{57} (2010),
121-134; Addendum to:
``Countably additive restrictions and limit
theorems in $(\ell)$-groups'', \textit{ibidem} \textbf{58}
(2011), 3-10.
\bibitem{bdpbjtm} A. Boccuto, X. Dimitriou and
N. Papanastassiou,
Brooks-Jewett-type theorems for the
pointwise ideal con\-ver\-gen\-ce  of measures
with values in $(\ell)$-groups, \textit{Tatra
Mt. Math. Publ.} \textbf{49} (2011), 17-26.
\bibitem{bdpschurfilters} A. Boccuto,
X. Dimitriou and N. Papanastassiou,
Schur lemma and limit theorems in lattice
groups with respect to filters, \textit{Math.
Slovaca} \textbf{62} (6) (2012), 1145-1166.
\bibitem{BRV} A. Boccuto, B. Rie\v{c}an and M.
Vr\'{a}belov\'{a},
\textit{Kurzweil-Henstock Integral in Riesz
Spaces,} Bentham Science Publ., U. A. E., 2009.
\bibitem{BSROMA} A. Boccuto and A. R.
Sambucini, The monotone integral
with respect to Riesz space-valued capacities, 
\textit{Rend. Mat. (Roma)} {\bf 16} (1996),
491-524.
\bibitem{CANDELOROPALERMO} D. Candeloro,
On the Vitali-Hahn-Saks,
Dieudonn\'{e} and Nikod\'{y}m theorems
(Italian), \textit{Rend. Circ. Mat. Palermo
Ser. II} Suppl. \textbf{8} (1985), 439-445.
\bibitem{CANDELOROLETTA} D. Candeloro and
G. Letta, On Vitali - Hahn - Saks and
Dieudonn\'{e} theorems (Italian), \textit{Rend.
Accad. Naz. Sci. XL Mem. Mat.} \textbf{9} (1)
(1985), 203-213.
\bibitem{cs} D. Candeloro and A. R. Sambucini,
Filter convergence and decompositions for
vector lattice-valued measures, \textit{Mediterranean J.
Math.} \textbf{12} (2015), 621-637.  
\bibitem{dobrakov1} I. Dobrakov, On submeasures I. 
\textit{Dissertationes Math.} \textbf{112} (1974), 5-35.
\bibitem{DRE2} L. Drewnowski,
Equivalence of Brooks - Jewett,
Vitali - Hahn - Saks and Nikodym theorems,
\it Bull. Acad. Polon. Sci. Ser. Sci. Math. Astronom.
Phys. \bf 20 \rm (1972), 725-731.
\bibitem{FREMLIN} D. H. Fremlin, A direct proof
of the Matthes-Wright integral extension theorem,
\textit{J. London Math. Soc.} \textbf{11} (2)
(1975), 276-284.
\bibitem{guariglia1990} E. Guariglia,
$k$-triangular functions on an orthomodular
lattice and the Brooks-Jewett theorem,
\textit{Radovi Mat.} \textbf{6} (1990), 241-251.
\bibitem{papvhs} E. Pap, The Vitali-Hahn-Saks Theorems for
$k$-triangular set functions, \textit{Atti Sem. Mat. Fis. Univ.
Modena} \textbf{35} (1987), 21-32.
\bibitem{PAPBOOK} E. Pap,
\textit{Null-Additive Set Functions,}
Kluwer Acad. Publishers/Ister Science,
Bratislava, 1995.
\bibitem{RN} B. Rie\v{c}an and T. Neubrunn,
\textit{Integral, Measure and
Ordering,} Kluwer Acad. Publ./Ister Science,
Dordrecht/Bratislava, 1997.
\bibitem{SAEKI} S. Saeki, The Vitali-Hahn-Saks
theorem and measuroids, \textit{Proc. Amer.
Math. Soc.} \textbf{114} (3) (1992), 775-782. 
\bibitem{SCHACHERMAYER}  W. Schachermayer,
On some classical measure-theoretic theorems for
non-sigma-complete Boolean algebras, 
\textit{Dissertationes Math.} \textbf{214} (1982), 1-33.
\bibitem{ventriglianapoli} F. Ventriglia,
Cafiero theorem for $k$-triangular functions
on an orthomodular lattice,
\textit{Rend. Accad. Sci. Fis. Mat. Napoli}
\textbf{75} (2008), 45-52. 
\bibitem{WANGKLIR} Z. Wang and G. J. Klir,
\textit{Generalized Measure theory,} Springer,
Berlin-Heidelberg-New York, 2009.
\end{thebibliography}
\end{document}